%% file: paper.tex
\documentclass[12pt]{amsart}
\font\emailfont=cmtt10

\usepackage{amsmath,amsthm,amsfonts,amscd,flafter,epsf,graphicx}
\usepackage[all,knot,arc]{xy}

\input{bas}

\title[{On invariants for Legendrian knots}]
{On invariants for Legendrian knots}

\author[Andr\'as I. Stipsicz]{Andr\'as I. Stipsicz}
\address{R\'enyi Institute of Mathematics, Budapest, Hungary and \newline
Department of Mathematics, Columbia University, New York 10027, USA
\newline
\indent{\emailfont{stipsicz@math.columbia.edu, stipsicz@renyi.hu}}}
\thanks{}

\author[Vera V\'ertesi]{Vera V\'ertesi}
\address{Institute of Mathematics , E\"otv\"os Lor\'and University,
Budapest, Hungary \newline
\indent{\emailfont{wera@szit.bme.hu}}}
\thanks{}

\begin{document}
\begin{abstract}
Suppose that $L$ is a null--homologous Legendrian knot in the contact
3--manifold $(Y, \xi )$. We determine the connection between the
sutured invariant $\EH (L)=\EH(Y-\nu (L),\xi|_{Y-\nu (L)})$ of $L$ and
the Legendrian invariant $\widehat{\mathcal{L}}(L)$ defined in
\cite{LOSS}. In particular, we derive a vanishing theorem for
$\LegInva (L)$ in the presence of Giroux torsion in the complement of
the knot, and reprove several known properties of the Legendrian
invariant from this perspective.
\end{abstract}
\maketitle

\section{Introduction}
\label{sec:intro}

A knot $L$ in a closed, contact 3--manifold $(Y, \xi )$ is \emph{Legendrian}
if the tangent vectors of the knot are contained by the contact 2--plane field
$\xi$. The knot $T$ is \emph{transverse}, if the (nonzero) tangent vectors are
not contained by $\xi$.  Legendrian and transverse knot theory has been shaped
by advances in convex surface theory \cite{EH} (showing that different looking
objects are actually equivalent) and by the introduction of various invariants
of these knots --- proving that different looking objects are, in fact,
different. Examples of such invariants are provided by Chekanov's differential
graded algebras and contact homology \cite{Ch,Ell}. More recently, Heegaard
Floer homology provided various sets of invariants: for knots in the standard
contact 3--sphere the combinatorial construction of knot Floer homology
through grid diagrams \cite{NOT, OSzT}, for null--homologous knots in general
contact 3--manifolds the Legendrian invariant of \cite{LOSS} and for general
Legendrian knots the sutured invariant of the knot complement \cite{HKM2}.

The aim of this paper is to set up a relation between these last two
invariants. To set the stage, recall that the Legendrian invariant $
\LegInva (L)$ of the null--homologous Legendrian knot $L\subset (Y,
\xi )$ defined in \cite{LOSS} takes its value in the knot Floer
homology group $\HFKa (-Y, L)$. (The theory admits a version where the
invariants are in the more refined group ${\rm {HFK}}^-(-Y, L)$, but
since the corresponding sutured theory is not developed yet, we will
deal only with the $\HFKa$--version in this paper.) In turn, the
sutured invariant $\EH (L)$ is defined as follows: consider the
Legendrian knot $L\subset (Y, \xi )$, and delete a standard
neighbourhood $\nu (L)$ of $L$ with convex boundary. The resulting
contact 3--manifold $Y-\nu (L)$ with convex boundary naturally admits
a balanced sutured 3--manifold structure $(Y-\nu (L), \Gamma )$, and hence by
\cite{Ju} it admits a sutured Floer homology $\SFH (Y-\nu(L), \Gamma
)$. According to \cite{HKM2} the contact structure on $Y-\nu (L)$
specifies an element $\EH (L)\in \SFH (-(Y-\nu(L)), -\Gamma)$, which we will
call the \emph{sutured invariant} of $L$. A relation between sutured
Floer homology and knot Floer homology obviously follows from their
definitions: suppose that $(Y-\nu (L), \Gamma )$ is the sutured
3--manifold with toric boundary we get by deleting a neighbourhood of
the (not necessarily Legendrian) knot $L$ and $\Gamma$ has two
(parallel) components. Then there is an obvious isomorphism
\[
\Psi \colon \SFH (Y-\nu(L), \Gamma )\to \HFKa (Y_{\Gamma}, L')
\]
where $Y_{\Gamma}$ is the Dehn filling of $Y-\nu (L)$ (and $L'$ is the
core of the Dehn filling) with slope given by the sutures $\Gamma$. In
general, $Y_{\Gamma}$ differs from $Y$ (and therefore $L'$ differs
from $L$). By attaching a specific contact $T^2\times [0,1]$ (a
\emph{basic slice}) to $Y-\nu (L)$, the composition of the map
\[
\Phi \colon \SFH (-(Y-\nu (L)),- \Gamma ) \to \SFH (-(Y-\nu (L)), \Gamma ')
\]
of \cite{HKM3} induced by this attachment and the above map $\Psi$
(applied to the suture $\Gamma '$ with components isotopic to the
meridian of the knot) gives a map
\[
F\colon \SFH (-(Y-\nu (L)), -\Gamma )\to \HFKa (-Y, L)
\]
for which we show the following:
\begin{theorem}\label{t:main}
Fix an orientation on the Legendrian knot $L$ and consider one of the
basic slices with boundary slopes given by the dividing set of
$\partial (Y-\nu (L))$ on $T^2\times \{ 0 \}$ and by the meridian of
$L$ on $T^2\times \{ 1\}$.  Then the map $F$ defined above maps $\EH
(L)$ to $\LegInva (L)$.
\end{theorem}
A more precise formulation of the theorem will be given in
Section~\ref{sec:conn} after basic slices and orientations have been
discussed.  A straightforward consequence of the above relation is the
following
\begin{corollary} \label{cor:girouxtorsion}
If the complement of a null--homologous Legendrian knot has
positive Giroux torsion then $\LegInva (L)$ vanishes.
\end{corollary}
\begin{remark}
{\rm The same corollary has been found recently by
D.~S.~Vela--Vick \cite{VV} using slightly different arguments.}
\end{remark}

To put this result in perspective, we recall that a knot type in the standard
contact 3--sphere is called \emph{Legendrian simple} if two Legendrian knots
of the given knot type and identical Thurston--Bennequin and rotation numbers
(for definitions of these invariants see \cite{Et}) are Legendrian isotopic.
The same notion generalizes to an arbitrary ambient contact 3--manifold $(Y,
\xi )$, with a caveat in the case when $\xi $ is overtwisted: in that case
Legendrian knots fall into two categories, depending on whether the knot
complement is overtwisted (in which case the knot is called \emph{loose}) or
--- although $\xi$ is overtwisted --- the knot complement is tight (in which
case the knot is \emph{non--loose} or \emph{exceptional}, cf.  \cite{EliFr}).
Obviously a loose and a non--loose knot cannot be isotopic, hence in
overtwisted contact 3--manifolds besides the equality of the
Thurston--Bennequin and rotation numbers we also require the equality of the
looseness of the two knots in defining simplicity. Non--simple non--loose
knots in a variety of overtwisted contact structures have been found in
\cite{LOSS}. There is, however, a simple way of constructing non--simple
non--loose knots \cite{Etn}: suppose that the knot complement contains an
incompressible torus (e.g., the knot type is a satellite in $S^3$) and
introduce Giroux torsion along the torus. Since this procedure does not change
the homotopy type of the 2--plane field, and $\xi$ is overtwisted by
assumption (and overtwisted structures are classified by their homotopy type),
after a suitable choice of the knot and the torus we get a Legendrian knot in
the same contact 3--manifold with different tight complement. (The
verification that the complement remains tight, and that the implementations
of different Giroux torsions result in different structures requires delicate
arguments \cite{Etn}.)  This method, in fact, can produce infinitely many
different Legendrian non--loose knots with the same numerical invariants in
these knot types \cite{Etn}. We say that $L\subset (Y, \xi )$ is
\emph{strongly non--loose} if $\xi$ is overtwisted and the knot complement is
tight with vanishing Giroux torsion. The knot type is \emph{strongly
  non--simple} if there are two strongly non--loose, smoothly isotopic knots
with equal numerical invariants which are not Legendrian isotopic. The same
simplicity/non--simplicity definition (with the strong adjective) carries
through verbatim for transverse knots (where the role of the numerical
invariants is played by the self--linking number of the transverse knot).  In
this sense, the result of \cite{LOSS} translates to
\begin{corollary}\label{cor:14}
The knot types of \cite[Theorem~1.7 and Corollary~1.8]{LOSS} are
strongly non--simple.
\end{corollary}
\begin{proof}
  The distinction of the Legendrian knots $L_i$ in \cite{LOSS} went by
  determining the Legendrian invariants $\LegInva (L_i)$, and since both were
  nonzero, Corollary~\ref{cor:girouxtorsion} implies that the knots $L_i$ are
  strongly non--loose, concluding the proof.
\end{proof}

Notice that in \cite{OSzT} the combinatorial theory provided two invariants of
$L$ (denoted by ${\widehat {\lambda}} ^{\pm }(L)$), while in \cite{LOSS} the
invariant $\LegInva (L)$ depended on an orientation of $L$ --- therefore an
unoriented Legendrian knot admitted two invariants $\LegInva (L)$ and
$\LegInva (-L)$ after an arbitrary orientation of $L$ was fixed. On the other
hand, the sutured theory provides a unique element for $L$. The discrepancy is
resolved by the observation that the map on sutured Floer homology induced by
the basic slice attachment is well--defined only up to a choice: with the
given boundary slopes there are two basic slices, and using one transforms
$\EH (L)$ into $\LegInva (L)$, while with the other choice the result will be
$\LegInva (-L)$ (after an orientation on $L$ is fixed).  In order to clarify
signs, we reprove a special case of \cite[Theorem~7.2]{LOSS} (only in the
$\HFKa$--theory) regarding the effect of stabilization of $L$ on $\LegInva$
and show

\begin{theorem}\label{thm:stabilization}
  Let $L$ be an oriented null--homologous Legendrian knot.  If $L^-$ (and
  $L^+$) denotes its negative (resp. positive) stabilization, then
  $\widehat{\mathcal{L}}(L^-)=\widehat{\mathcal{L}}(L)$ and
  $\widehat{\mathcal{L}}(L^+)=0$.
\end{theorem}

Notice that the invariance of $\LegInva$ under negative stabilization means
that, in fact, it is an invariant of the transverse isotopy class of the
positive transverse push--off of the Legendrian knot $L$.  By this definition
the extensions of Corollaries~\ref{cor:girouxtorsion} and \ref{cor:14} to the
transverse case are easy exercises.  For further results regarding transverse
knots using these invariants see \cite{NOT, OS}. In fact, in \cite{OS} the
distinction of various Legendrian and transverse Eliashberg--Chekanov (aka
twist) knots and 2--bridge knots was carried out by computing their
$\LegInva$--invariants. As a corollary, Theorem~\ref{t:main} readily implies

\begin{corollary}
The complement of the Eliashberg--Chekanov knot $E_n$ (which is the
2--bridge knot of type $\frac{2n+1}{2}$) for odd $n$ admits at least
$\lceil \frac{n}{4} \rceil$ different tight contact structures
(distinguished by the sutured invariant) with convex boundary and
dividing set $\Gamma$ of two components with slope 1.  \qed
\end{corollary}

Performing contact $(-1)$--surgery along a Legendrian knot $L$ gives a
well--defined contact structure $\xi_{-1}$ on the surgered 3--manifold
$Y_{-1}$. The core $L'$ of the glued--back solid torus is a Legendrian knot in
$(Y_{-1}, \xi _{-1})$. Suppose that $L^\prime$ is null--homologous in
$Y_{-1}$.  Using the sutured invariant we deduce
\begin{theorem} \label{thm:surgery}
Under the circumstance described above $\LegInva(L)\neq 0$ implies
$\LegInva(L^\prime)\neq 0$.
\end{theorem}

The paper is organized as follows. In Section \ref{sec:invariant} we review
the basic definitions we need about contact structures. Section~\ref{sec:HF}
gives a short description of sutured Heegaard Floer homology and the
definition of the Legendrian invariants.  In Section~\ref{sec:conn} we state a
precise version of Theorem \ref{t:main} and prove it together with the
consequences given in the Introduction.

\smallskip

\noindent {\bf Acknowledgements:} We would like to thank Ko Honda and
Paolo Ghiggini for helpful discussions. AS acknowledges support from
the Clay Mathematics Institute. AS was also partially supported by
OTKA 49449 and by Marie Curie TOK project BudAlgGeo. VV was supported
by NSF grant number FRG-0244663 and OTKA 49449 and 67867. VV was also
supported by ``Magyar \'Allami E\"otv\"os \"Oszt\"ond\'\i j''.

\section{Contact preliminaries}\label{sec:invariant}

\subsection{Contact 3-Manifolds}\label{subsec:contact}
A surface $\Sigma$ in the contact 3--manifold $(Y, \xi )$ is
\emph{convex} if there is a contact vector field $X$ defined near
$\Sigma$ which is transverse to $\Sigma$. The set of points $p\in
\Sigma$ where $X_p\in \xi _p$ is usually denoted by $\Gamma$ and
called the \emph{dividing set} of the convex surface $\Sigma$. It
turns out that $\Gamma$ is an embedded 1--manifold, partitioning
$\Sigma$ into $\Sigma _+$ and $\Sigma _-$, and the contact structure
$\xi$ is determined by $\Gamma$ near $\Sigma$. For a more complete
treatment of the subject, see \cite{Et, EH, OzS}.

\begin{figure}[ht]
\centering
\includegraphics[scale=0.5]{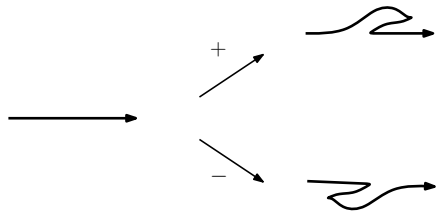}
\caption{Positive and negative stabilization.}
\label{fig:stabi}
\end{figure}

Suppose that $L$ is an oriented null--homologous Legendrian knot in
the contact 3--manifold $(Y, \xi )$. Let $S$ be a Seifert surface of
$L$ in convex position. Orient $S$ such that its boundary orientation
gives the orientation for $L$. The rotation number then can be
computed as $\rot(L)=\chi(S_+)-\chi(S_-)$. Define the \emph{negative}
and \emph{positive stabilizations} $L^-$ and $L^+$ by modifying $L$
near a point as it is depicted by Figure~\ref{fig:stabi}.  The effect
of a positive (resp. negative) stabilization on the numerical
invariants of $L$ can be easily computed as
\[
\tb (L^{\pm })=\tb (L)-1 \qquad {\mbox{ and }} \qquad \rot (L^{\pm })=
\rot (L)\pm 1.
\]
Notice that the sign of the stabilization makes sense only after
fixing an orientation for the Legendrian knot.

\subsection{Sutured 3--manifolds}
A \emph{sutured 3--manifold} is a pair $(Y, \gamma )$ where $Y$ is a
compact, oriented 3--manifold with boundary and $\gamma \subset
\partial Y$ is a disjoint set of embedded tori and annuli.  Every component of
$R (\gamma )=\partial Y - \gamma$ is oriented, and $R_+$ (resp.  $R_-$) is the
union of those components where the normal vector points out (resp. in) $Y$.
The sutured manifold is called \emph{balanced} if all sutures are annular, $Y$
has no closed components, every boundary component admits a suture and $\chi
(R_+)=\chi (R_-)$ on every component of $Y$. As is customary, annular sutures
are symbolized by the homologically nontrivial simple closed curves they
contain, the collection of which is denoted by $\Gamma$.  Without confusion,
the term ``suture'' will also refer to these curves, and sometimes to their
union $\Gamma$. The suture $\Gamma$ is oriented as the boundary of $R_+\subset
\partial Y$. We will consider only balanced sutured manifolds in this paper.

\subsection{Partial Open Books}\label{subsec:partialopenbook}
Partial open books are generalizations of open books for 3--manifolds
with boundary. This notion was introduced by Honda, Kazez and Mati\'c
in \cite{HKM2}, see also \cite{EO1,EO2}.

\begin{definition} {\rm An \emph{abstract partial open book} is a triple $(S,
    P, h)$ where $S$ is a connected surface with boundary, $P$ is a proper
    subsurface of $S$ which is a union of 1--handles attached to $S-P$, and $h
    \colon P\to S$ is an embedding that restricts to the identity near the
    boundary $\partial P\cap \partial S$.}
\end{definition}

\begin{figure}[ht]
\centering
\includegraphics[scale=0.5]{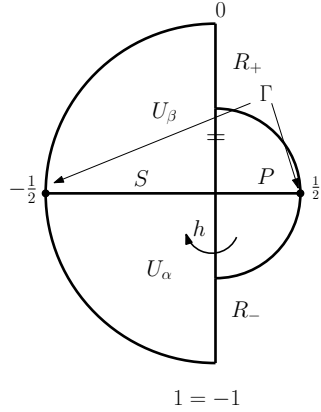}
\caption{Schematic picture of a partial open book decomposition.}
\label{fig:partialopenbook}
\end{figure}

A partial open book defines a 3--manifold $Y$ with boundary as follows. First
construct the handlebody $S\times[-1,0]/\sim$ and the compression--body
$P\times[0,1]/\sim$, where $(x,t)\sim(x,t^\prime)$ for $x\in \partial S$ and
$t,t^\prime\in [-1,1]$. (Note that on $P\times[0,1]$ we just contract the
points with first coordinate in $\partial P\cap\partial S$.)  Then glue them
together with the maps $P\times\{0\}\hookrightarrow S\times\{0\}$ and $h\colon
P\times\{1\} \to S\times\{-1\}$.  A schematic picture of $Y$ is given by
Figure~\ref{fig:partialopenbook}. The resulting 3--manifold naturally carries
the structure of a balanced sutured manifold: take $\Gamma =
\overline{\partial S - \partial P}\times \{-\frac12\}\cup -\overline{(\partial
  P -\partial S)}\times \{\frac12\}\subset \partial Y$.  Now $R_+=\overline{S-
  P}\times\{0\}$, $R_-=\overline{S- h(P)}\times\{-1\}$, consequently $\chi
(R_+)= \chi (R_-)$ follows at once.

Both the handlebody $S\times[-1,0]/\sim$ and the compression--body $P\times
[0,1]/\sim$ admit unique tight contact structures with convex boundary and
dividing set $\partial S$ (and $\partial P$, resp), cf.  \cite{EO2, Tor}. As
the dividing sets match up, we can glue these contact structures to obtain a
contact structure $\xi$ on $Y$ with dividing set $\Gamma$ on the convex
boundary $\partial Y$.  In this sense a partial open book decomposition
determines a contact structure with convex boundary (inducing the dividing set
given by the sutures associated to the partial open book).

The partial open book decomposition naturally induces a Heegaard
decomposition of $Y$ with the compression bodies
$U_\alpha=P\times[\frac12,1]\cup S\times[-1,-\frac12]$ and $U_\beta=
S\times[-\frac12,0]\cup P\times[0,\frac12]$, divided by the Heegaard
surface $\Sigma=\partial U_\alpha= S\times\{-\frac12\}\cup
-P\times\{\frac12\}$. Consistently with the sutured 3--manifold
structure, the boundary of $U_\alpha$ (and $U_\beta$, resp.) consists of the union
of $\Sigma$ (resp. $-\Sigma$), $R_-$ (resp. $R_+$) and a collar
neighbourhood for $\Gamma$; furthermore $\Gamma=\partial
\Sigma(=\partial R_+=-\partial R_-)$.

Every contact 3--manifold with convex boundary $(Y,\xi)$ admits a partial open
book decomposition that is compatible with $\xi$ in the above sense, cf.
\cite{HKM2}. To see this, consider a contact cell--decomposition for $Y$ whose
1--skeleton $C$ is a direct product near the boundary $\partial Y$ and
intersects the boundary on the dividing set. As Legendrian arcs have standard
neighbourhood, there is a neighbourhood $\nu (C)$ of $C$ with convex boundary
and with dividing curves of two components. The dividing curve separates
$-\partial \nu (C)$ into a positive and a negative part $(-\partial \nu
(C))_+$ and $(-\partial \nu (C))_-$. Setting $P=(-\partial \nu (C))_+$ the
neighbourhood $\nu (C)$ can be written as $P\times [0,1]/\sim$. As $C$ was the
1--skeleton of a contact cell--decomposition, $Y-\nu (C)$ is product
disk--decomposable: it is divided by the 2--cells of the contact
cell--decomposition (that are disks with $\mathrm{tb}=-1$) to a union of tight
contact 3--balls. Thus for $S=\partial (Y-\nu (C))_+$ the handlebody $Y-\nu
(C)$ can be written as $Y-\nu (C)=S\times[-1,0]/\sim$, and
$P=(-\partial(\nu(C)))_+\subseteq (\partial(Y-\nu(C)))_+=S$. Note that by
construction $\xi|_{Y-\nu (C)}$ is tight, its boundary $\partial(Y-\nu (C))$
is convex, and the dividing set $\Gamma_{\partial(Y-\nu (C))}$ is isotopic to
$\partial S\times\{0\}$.

\subsection{Bypass attachment} \label{subsec:bypass}
Next we review the change of the partial open book decomposition after
a bypass is attached along a Legendrian curve $c$ on the boundary. For
a complete discussion of
bypass attachments see \cite{Honda}. The
considerations below already appeared in \cite[Example 5]{HKM2}.

\begin{figure}[ht]
\centering
\includegraphics[scale=0.7]{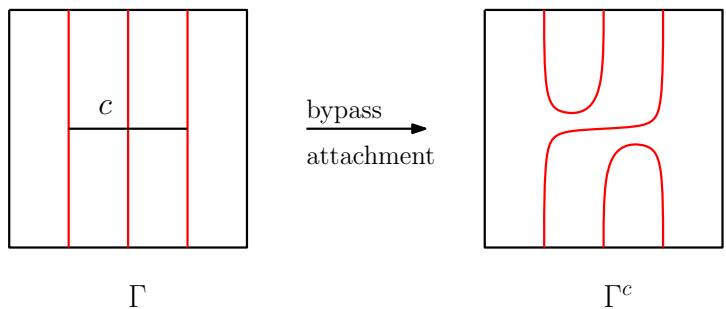}
\caption{Bypass attachment.}
\label{fig:bypass1}
\end{figure}

Let $(Y,\partial Y,\xi)$ be a contact 3--manifold with convex
boundary. Suppose that we are given a Legendrian arc $c\subset
\partial Y$ that starts and ends on the dividing set $\Gamma_{\partial
Y}$ and intersects $\Gamma_{\partial Y}$ in one additional
point. Attaching a \emph{bypass} along $c$ is --- roughly speaking --- the
attachment of the neighbourhood of a ``half overtwisted disk''. This is a
disk $D$ with boundary $\partial D=c\cup d$, where $\partial D\cap
\partial Y=c$, and the dividing curve on $D$ consists of a single arc with
both of its endpoints on $c$. The resulting manifold is diffeomorphic to $Y$
with contact structure $\xi^c$, and the dividing curve $\Gamma$ is changed in
the neighbourhood of $c$ to $\Gamma^c$ as it is shown on Figure
\ref{fig:bypass1}.

\begin{figure}[ht]
\centering
\includegraphics[scale=0.7]{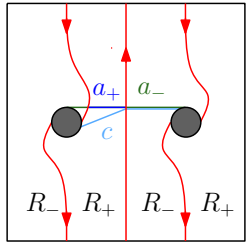}
\caption{The grey areas indicate the attaching regions of the 1--handle. The
  attaching curve for the 2--handle is $a=a_+\cup a_-$ and $a_-$ is assumed to
  go parallel to the core of the 1--handle in the negative region.}
\label{fig:bypass}
\end{figure}

Take a partial open book decomposition for $(Y,\partial Y,\xi)$ coming from a
contact cell--decomposition whose 1--skeleton $C$ misses the attaching arc
$c$.  Let $c_\pm=c\cap R_\pm$. Under the identification of $Y-\nu (C)$ with
$S\times[-1,0]/\sim$ both $c_+$ and $c_-$ are arcs on $S$.  The bypass
attachment can be thought of consisting of a 1--handle attachment with core
$d$ followed by a canceling 2--handle attached along the curve $a=a_+\cup a_-$
of Figure~\ref{fig:bypass}. The contact cell--decomposition can be extended to
the new manifold $(Y^{\prime}, \xi ')$ (where $Y'$ is, in fact, diffeomorphic
to $Y$) by including the cocore of the 2--handle in the 1--skeleton. Thus
$C^{\prime}=C\cup(\textrm{cocore of the 2--handle})$ and the page $S^{\prime}$
of the partial open book decomposition resulting from this contact
cell--decomposition will be equal to $S\cup (-\partial \nu
(d))_+=S\cup(\textrm{1--handle})$.  Denote the intersection of the attaching
circle of the canceling 2--handle with the positive and negative parts of
$\partial(Y\cup \{\textrm{1--handle}\})$ by $a_\pm=a\cap R^{\prime}_\pm$. As
it is depicted in Figure~\ref{fig:bypass}, the arc $a_+$ can be pushed off to
lie entirely in the boundary of the old manifold $Y$, thus $a_+\subset R_+$.
Note that $c_+$ and $a_+$ are isotopic. They have one endpoint that agrees
with the endpoint of both $c_-$ and $a_-$, and the other one is moved in the
direction given on $\Gamma_{\partial Y}$ as the boundary of $R_+$. These
curves can again be thought of as being on $S$. Now $R_+^{\prime}=(R_+-\nu
(a_+))\cup (\partial (-\nu (d)))_+$, thus $P^\prime=P\cup \nu (a_+)$. The
mondoromy $h^{\prime}$ remains the same on $P$, so we only need to understand
it on $a_+$. To push $a_+$ through $\nu (C)$ we just have to push it through
the newly attached 1--handle, so $h^\prime(a_+)=a_-$. The arc $a_-$ can be
split to two subarcs $a_-\cap S$ and the core of the 1--handle in $S^\prime$.

\subsection{Basic slices}
We give a short description of basic slices defined by Honda \cite{Honda}.
Suppose that $\xi$ is a contact structure on $T^2\times[0,1]$ with convex
boundary with two--component dividing curves on each of its boundary
components.  The dividing curves are homotopically nontrivial and parallel.
Fix a trivialization for $T^2$ as $\bfr ^2/\bfz ^2$ and let $s_i$ denote the
slope of the dividing curves on $T^2\times\{i\}\ $ ($i\in\{0,1\}$).  The
contact 3--manifold $(T^2\times [0,1],\xi)$ is called \emph{minimally
  twisting} if every convex torus parallel to the boundary has slope $s$ in
$[s_1,s_0]$. (By $[s_1,s_0]$ we mean $[s_1,\infty]\cup[-\infty,s_0]$ if
$s_1\ge s_0$.)  A \emph{basic slice} is a minimally twisting tight contact
structure $(T^2\times[0,1],\xi)$, with convex boundary and with two dividing
curves on each $T^2\times\{i\}$ and boundary slopes $s_0$ and $s_1$ forming an
integral basis for $\mathbb{Z}^2$.  For fixed boundary conditions (up to
isotopy) there are two basic slices distinguished by their relative Euler
class, which differ by their sign; there is no canonical positive or negative
choice.

One way to obtain a basic slice is by gluing a bypass to an
$I$--invariant neighbourhood of a convex $T^2$ with two dividing
curves. For a given slope of the attaching curve there are two ways of
attaching a bypass corresponding to the two different basic slices,
cf. Figure~ \ref{fig:meridian}. Any basic slice can be obtained by
this construction.

Suppose that $(T^2\times[0,1],\xi_0)$ and $(T^2\times[1,2],\xi_1)$ are basic
slices with boundary slopes $s_i$ on $T^2\times\{i\}\ $ ($i\in\{0,1,2\}$). As
the dividing curves match up on $T^2\times\{1\}$, we can glue them together to
obtain $(T^2\times[0,2],\xi =\xi_0\cup\xi_1)$. If in addition we require that
the shortest representatives of $s_0$ and $s_2$ give an integral basis for
$\mathbb{Z}^2$ and $[s_0,s_1]\cup[s_1,s_2]\neq [-\infty,\infty]$, then
$(T^2\times[0,2],\xi)$ is minimally twisting. It is either overtwisted or a
single basic slice depending on whether the basic slices
$(T^2\times[0,1],\xi_0)$ and $(T^2\times[1,2],\xi_1)$ have the same or
opposite signs. Note that ``having the same sign'' makes sense in this
setting, once we require the trivialization of $\xi_0$ and $\xi_1$ to agree
over $T^2\times\{1\}$.

\section{Heegaard Floer invariants}\label{sec:HF}
In \cite{OSz1, OSz2} invariants of closed, oriented 3--manifolds have been
introduced.  In the simplest version, these invariants are given as follows.
Suppose that the 3--manifold $Y$ is given by a Heegaard diagram $(\Sigma _g,
\alphas , \betas )$, where $\Sigma _g$ is a genus--$g$ surface, the $g$
$\alpha$--curves $\alphas = \{ \alpha _1, \ldots , \alpha _g\}$ correspond to
belt circles of 1--handles, while the $g$ $\beta$--curves $\betas = \{ \beta
_1 , \ldots , \beta _g \}$ to attaching circles of 2--handles in a handle
decomposition of $Y$ with a unique 0-- and 3--handle. In particular, the
$\alpha$-- (and similarly the $\beta$--) circles are disjoint, and linearly
independent in homology. By fixing a base point $w\in \Sigma _g$ in the
complement of all the $\alpha$-- and $\beta$--curves, the chain complex $(\CFa
(Y),
\partial )$ is defined as follows: consider the $\bfz _2$--vector
space $\CFa (Y)$ freely generated by the intersections $\Ta \cap \Tb
\subset {\rm Sym} ^g (\Sigma _g )$, where the tori $\Ta $ and $\Tb$
are the products of the $\alpha$-- and $\beta$--curves,
respectively. The boundary operator $\partial $ is defined by counting
holomorphic disks in ${\rm Sym} ^g (\Sigma _g)$ (for an appropriate
choice of almost complex structure) connecting intersection points of
$\Ta$ and $\Tb$ which avoid the divisor $V_w =\{ w \}\times {\rm Sym}
^{g-1}(\Sigma _g )$. If $(\Sigma _g , \alphas , \betas )$ satisfy the
technical condition of admissibility (which can always be arranged by
suitable isotopies, cf. \cite{OSz2}) then the homology $\HFa (Y)$ of
the resulting chain complex is a diffeomorphism invariant of $Y$.

Variants of this construction provide invariants for knots and for sutured
3--manifolds, as will be outlined below. First, the choice of another point
$z\in \Sigma _g$ in the complement of the $\alpha$-- and the $\beta$--curves
determines a knot $K\subset Y$, and by taking $\CFKa (Y , K )= \CFa (Y)$ and
modifying $\partial $ to $\partial _K$ by only allowing holomorphic disks
avoiding both $V_w$ and $V_z$ we get a chain complex $(\CFKa (Y, K),
\partial _K)$, with homology the \emph{knot Floer homology} group
$\HFKa (Y, K)$. As it is shown in \cite{OSzknot, Ras}, for $K$
null--homologous in $Y$ this homology group will be an invariant of
the pair $(Y, K)$.

Suppose now that $\Sigma$ is a compact surface with nonempty boundary. Then by
fixing $k$ linearly independent (in homology) and disjoint $\alpha$-- (and
similar $\beta$--) circles, the attachment of the appropriate handles gives a
balanced sutured 3--manifold with sutures being equal to $\partial \Sigma$. In
fact, every balanced sutured 3--manifold arises in this way.  The previous
scheme applies verbatim (without even the choice of base points) and provides
a chain complex $(\SFC (Y, \Gamma ) , \partial _{\Gamma})$, ultimately
defining the \emph{sutured Floer homology} group $\SFH (Y, \Gamma )$, which
has been shown to be an invariant of the sutured 3--manifold \cite{Ju}.

If $\Sigma$ has exactly two boundary components and ${\overline {\Sigma}}$
denotes the capped--of closed surface, and if the number of attaching curves
$k$ equals to the genus of $\Sigma$ and the curves are homologically
independent in ${\overline {\Sigma }}$, then the corresponding sutured
3--manifold has toric boundary with a 2--component suture, and by placing two
marked points on the caps we get an identification
\[
\Psi\colon \SFH (Y , \Gamma )\to \HFKa (Y_{\Gamma}, L^\prime),
\]
where $Y_{\Gamma}$ is the result of Dehn filling of $Y$ with slope given by a
component of $\Gamma$ and $L^\prime$ is the core of the glued--up
solid torus.

\subsection*{The contact invariant}\label{subsec:contactinvariant}
Suppose that $(Y, \xi )$ is a contact 3--manifold with convex boundary, and
consider a partial open book compatible with $\xi$.  Let $\{b_1,\dots, b_k\}$
be a basis for $H_1(P,\partial S\cap \partial P)$. The disks swept out by the
$b_i$'s in the $U_\beta$ handlebody have boundaries
$\beta_i=b_i\times\{\frac12\}\cup b_i\times\{-\frac12\}$. Isotope each $b_i$
to an arc $a_i$ that intersects it transversely in a single point, and whose
endpoints are moved in the direction given by the boundary orientation of
$-P$. In the $U_\alpha$ handlebody $a_i$ sweeps out a disk with boundary
$\alpha_i=a_i\times\{\frac12\}\cup h(a_i)\times \{-\frac12\}$, providing a
Heegaard diagram $(\Sigma,\alphas, \betas)$ for $(Y,\Gamma)$. The single
intersection point $\y=(a_i\cap b_i)$ on $P\times\{\frac12\}$ can be shown to
represent a cycle in $\SFC(-\Sigma,\alphas,\betas)$, thus it defines an
element $\EH (Y,\xi)$ in $\SFH(-Y,-\Gamma)$. (Notice the reversal of
orientation of the Heegaard surface $\Sigma$.) As has been proven by Honda,
Kazez and Mati\'c \cite{HKM2}, this element is independent of the choices made
throughout its definition and gives the invariant $\EH (Y, \xi )$ of the
contact structure $(Y,\xi)$.  In the special case when the contact 3--manifold
with convex boundary is given as the complement of a standard neighbourhood of
a Legendrian knot in a closed contact 3--manifold $(Y, \xi )$, the resulting
element will be denoted by $\EH (L)$. Note that by the Legendrian
Neighbourhood Theorem, in this case $\Gamma$ consists of two parallel simple
closed curves in $\partial(Y-\nu(L))$.

\subsection*{The Legendrian invariant}\label{legendrianinvariant}
Consider an oriented, null--homologous Legendrian knot in the closed contact
3--manifold $(Y, \xi )$.  There is an open book decomposition of $Y$
compatible with $\xi$ containing $L$ on one of its pages
$S=S\times\{\frac12\}$. Consider a properly embedded arc $b_1$ in $S$
intersecting $L$ exactly once. The disk $b_1\times[0,1]$ is a meridional disk
for $L$. Orient $b_1$ so that the boundary orientation of
$\partial(b_1\times[0,1])=-b_1\times\{0\}\cup b_1\times\{1\}$ agrees with the
natural orientation of the meridian for $L$. (Such an oriented arc $b_1$ will
be called a \emph{half--meridian} of $L$.) With these conventions the
orientation of $S$ coincides with the orientation induced by $(b_1,L)$. Our
setup here will be slightly different from the one used in \cite{LOSS}, but
the resulting Heegaard diagram and the element specified in it will be
actually the same already on the chain--level.

\begin{figure}[ht]
\centering
\includegraphics[scale=0.7]{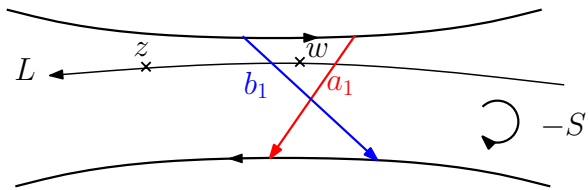}
\caption{The placement of the basepoints.}
\label{fig:page}
\end{figure}

Pick a basis $\{b_1,\dots, b_g\}$ of $H_1(S, \partial S)$ such that $b_1$ is a
half--meridian of $L$. Isotope all the $b_i$'s to $a_i$'s as before and place
the basepoint $z$ in the ``big'' region that is not swept out by the isotopies
of the $b_i$, and put $w$ between $b_1$ and $a_1$. This can be done in two
essentially different ways, and exactly one of them corresponds to the chosen
orientation of $L$. If $b_1$ is oriented as described above, $w$ should be
placed close to the tail of $b_1$, cf. Figure \ref{fig:page}.
The single intersection point $(a_i\cap b_i)$ on $S\times\{\frac12\} \subset
-\Sigma$ is an element in $\CFKh(-\Sigma ,\alphas ,\betas ,z,w)$ and the
choice of $z$ assures that it is a cycle, hence it defines an element
$\LegInva(L)$ in $\HFKh(-Y,L)$.  As it was shown in \cite{LOSS}, the homology
class $\LegInva(L)$ is an invariant of the oriented Legendrian knot $L\subset
(Y,\xi)$.

\section{Connection between the invariants}
\label{sec:conn}
Let $L$ be a Legendrian knot in a closed contact 3--manifold $(Y,\xi)$. The
two invariants $\EH (L)= \EH(Y-\nu (L),\xi|_{Y-\nu (L)})\in\SFH(-(Y-\nu (L)),
-\Gamma_{\partial(Y-\nu (L))})$ and $\LegInva(L)\in \HFKh(-Y,L)$ introduced
above lie in two different groups, but if we change the suture on
$\partial(Y-\nu (L))$ to two meridians $-m\cup m$ of $L$, the sutured Floer
homology $\SFH(-(Y-\nu (L)),-m \cup m)$ can be identified with $\HFKh(-Y,L)$.
This modification of the suture can be achieved by attaching a basic slice to
the sutured 3--manifold $Y-\nu (L)$, and according to \cite{HKM3} there is a
map corresponding to this attachment. More generally:

\begin{theorem}[Honda--Kazez--Mati\'c, \cite{HKM3}, cf. also
\cite{GH}]\label{thm:gluing}
Suppose $(Y^\prime, \Gamma^\prime)$ is a balanced sutured submanifold of the
balanced sutured 3-manifold $(Y,\Gamma)$ and all components of $Y -
\textrm{int}(Y^\prime)$ intersect $\partial Y$. Let $\xi$ be a contact
structure on $Y - \textrm{int}(Y^\prime)$ so that $\partial Y \cup
\partial Y^\prime$ is convex with respect to $\xi$ and with dividing
set $\Gamma \cup \Gamma^\prime$. Then there is a natural linear map
\[
\Phi_\xi \colon \SFH(-Y^\prime,-\Gamma^\prime) \to \SFH(-Y,-\Gamma),
\]
induced by $\xi$. Moreover, if $Y^\prime$ is endowed with the contact
structure $\xi^\prime$ such that
$\Gamma_{(Y^\prime,\xi^\prime)}=\Gamma^\prime$ then
\[
\Phi_\xi(\EH (Y^\prime,\xi^\prime))=\EH (Y,\xi^\prime\cup\xi).
\]
\qed
\end{theorem}

We will apply this theorem in the special case when $\partial Y'$ and
$\partial Y$ are both 2--tori, $Y-{\rm {int}}Y'=T^2 \times [0,1]$ and the
contact structure on the difference is a basic slice.  The dividing set is
given on $\partial (T^2\times [0,1])$ by the dividing set of $\partial Y$ (on
$T^2\times \{ 0 \}$) and by the meridians of $L$ (on $T^2\times \{1\}$); there
are two basic slices with the given boundary slopes. Notice that the
attachment of the basic slice is actually equivalent to the attachment of a
single bypass.

\begin{figure}[ht]
\centering \includegraphics[scale=0.7]{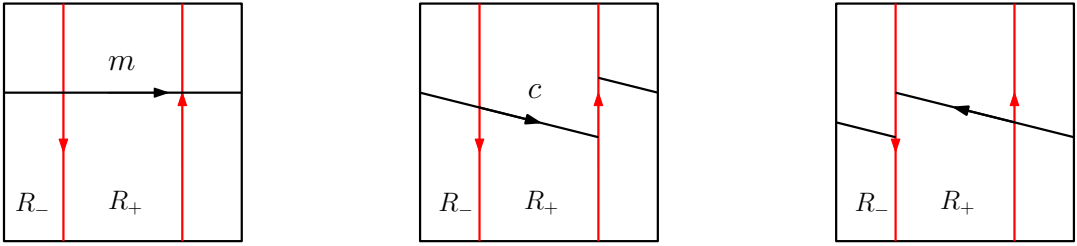}
\caption{Bypass attachments to obtain meridians.}
\label{fig:meridian}
\end{figure}

Trivialize $\partial (Y-\nu (L))$ with the meridian $m$ and the contact
framing $l$, hence the dividing curves have slope $\infty$.  The new dividing
curve after attaching a bypass along any arc with slope between $-1$ and $0$
has slope $0$.  Up to isotopy there are only two different attachments (of
opposite sign) depicted on Figure \ref{fig:meridian}; these are the two
different bypass attachments corresponding to the two different basic slices.
These attaching curves together with the arcs of the dividing curves form an
oriented curve on $\partial(Y-\nu (L))$, one of them represents $m$ the other
one represents $-m$. Denote the former one by $c$.
\begin{theorem}\label{thm:transf}
The map
\[
\Phi^c \colon \SFH(-(Y-\nu (L)),-\Gamma_{\partial(Y-\nu (L))})\to
\SFH(-(Y-\nu (L)),-m\cup m)
\]
induced by the basic slice attachment along $c$ maps $\EH (L)$ to the
class which is identified with $\LegInva (L)$ under the identification
\[
\Psi \colon \SFH(-(Y-\nu (L)),-m\cup m) \to \HFKh(-Y,L).
\]
\end{theorem}
\begin{proof}
  Let $(S,g)$ be an open book for $(Y,\xi)$ that contains $L$ homologically
  essentially on one of its pages. Set $P=S- \nu _S (L)$ (where $\nu _S (L)$
  denotes the tubular neighbourhood of $L$ in $S$) and $h=g|_{P}$. We claim
  that the partial open book $(S,P,h)$ describes $(Y- \nu (L),\xi|_{Y-\nu
    (L)})$. Indeed, topologically the 3--manifold corresponding to this
  abstract partial open book is $(S\times[-1,0]/\sim) \cup
  (P\times[0,1]/\sim)$, which is equal to
\[
(S\times[-1,1]/\sim)-(\nu _S (L)\times[0,1])=Y-\nu (L).
\]
The contact structure on $S\times[-1,0]/\sim$ is the same, while on
$P\times [0,1]/\sim$ (which is a subset of $S\times [0,1]/\sim$) it is
obviously tight. If we round the corners we get that the dividing
curve is $\Gamma_{\partial(Y-\nu (L))}$, so the dividing curve on
$P\times [0,1]/\sim$ must be $\partial P$.

Take a basis $\{b_1,\dots,b_k\}$ of $S$ subordinated to $L$, such that
$b_1$ is the half--meridian of $L$.  Then the left hand side of
Figure~\ref{fig:heegaardiagram} depicts the corresponding Heegaard
diagram $(-\Sigma,\{\alpha_1,\dots,\alpha_k\},
\{\beta_1,\dots,\beta_k\},w,z)$ for $(-Y,L)$. Here
$\Sigma=S\times\{\frac12\}\cup -S\times\{-\frac12\}$ and the
intersection point $\ensuremath{\x=(a_i\cap b_i)_{i=1}^{k}}$
represents the Legendrian invariant $\LegInva(L)$ in
$\HFKh(-Y,L)$. The basis for $H_1(P,\partial S\cap \partial P)$ is
$\{b_2,\dots,b_k\}$ while the Heegaard surface is
$-\Sigmah=P\times\{\frac12\}\cup -S\times\{-\frac12\}$. The
corresponding Heegaard diagram for $(-(Y-\nu
(L)),-\Gamma_{\partial(Y-\nu (L))})$ is
$(-\Sigmah,\{\alpha_2,\dots,\alpha_k\}, \{\beta_2,\dots,\beta_k\})$
which is depicted on the right hand side of
Figure~\ref{fig:heegaardiagram}. By definition
$\ensuremath{\y=(a_i\cap b_i)_{i=2}^{k}}$ represents the contact
invariant $\EH (L)\in \SFH(-(Y-\nu (L)),-\Gamma_{\partial(Y-\nu
(L))})$.

\begin{figure}[ht]
\centering
$\begin{array}{c@{\hspace{1cm}}c}
\includegraphics[scale=0.35]{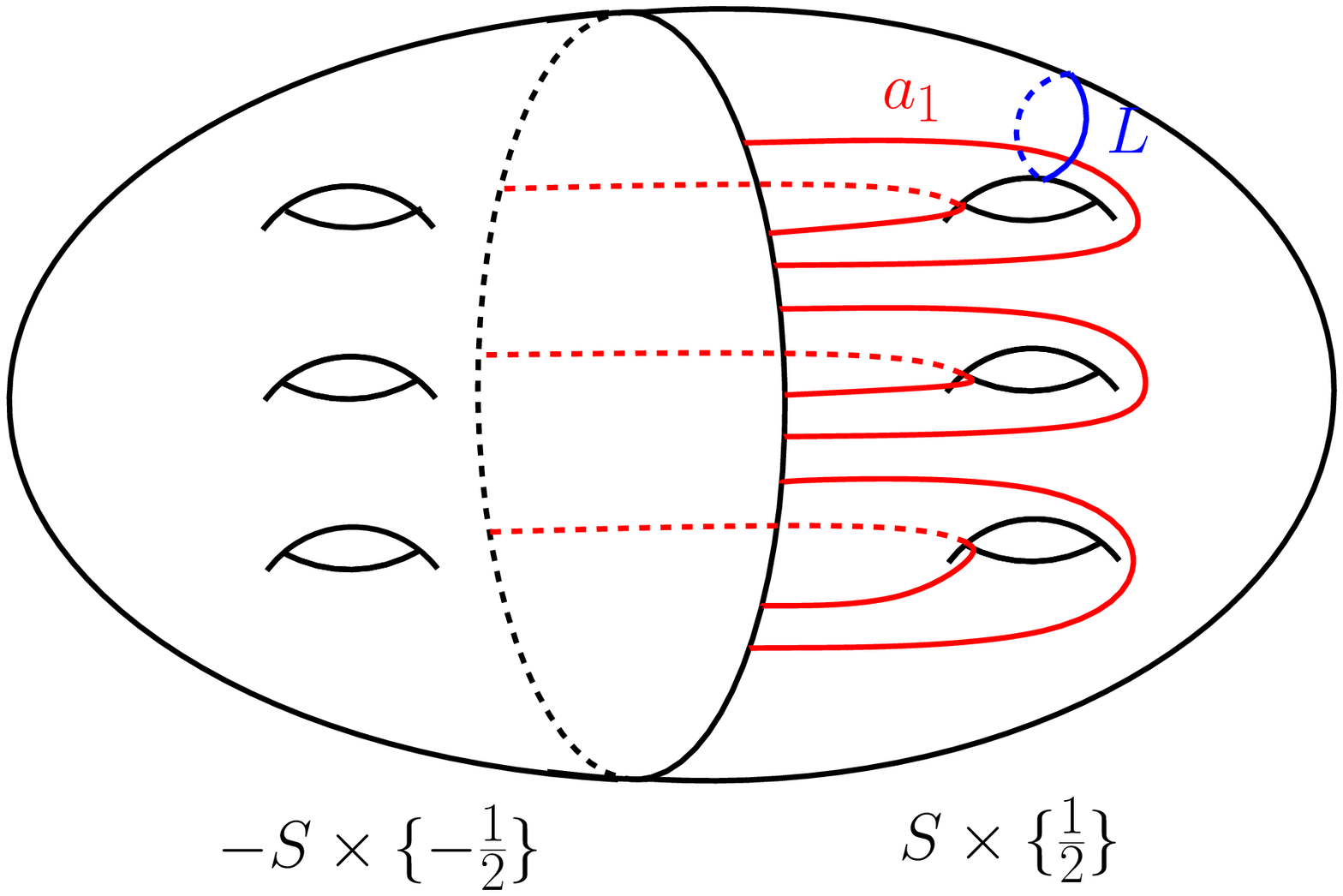}
&
\includegraphics[scale=0.35]{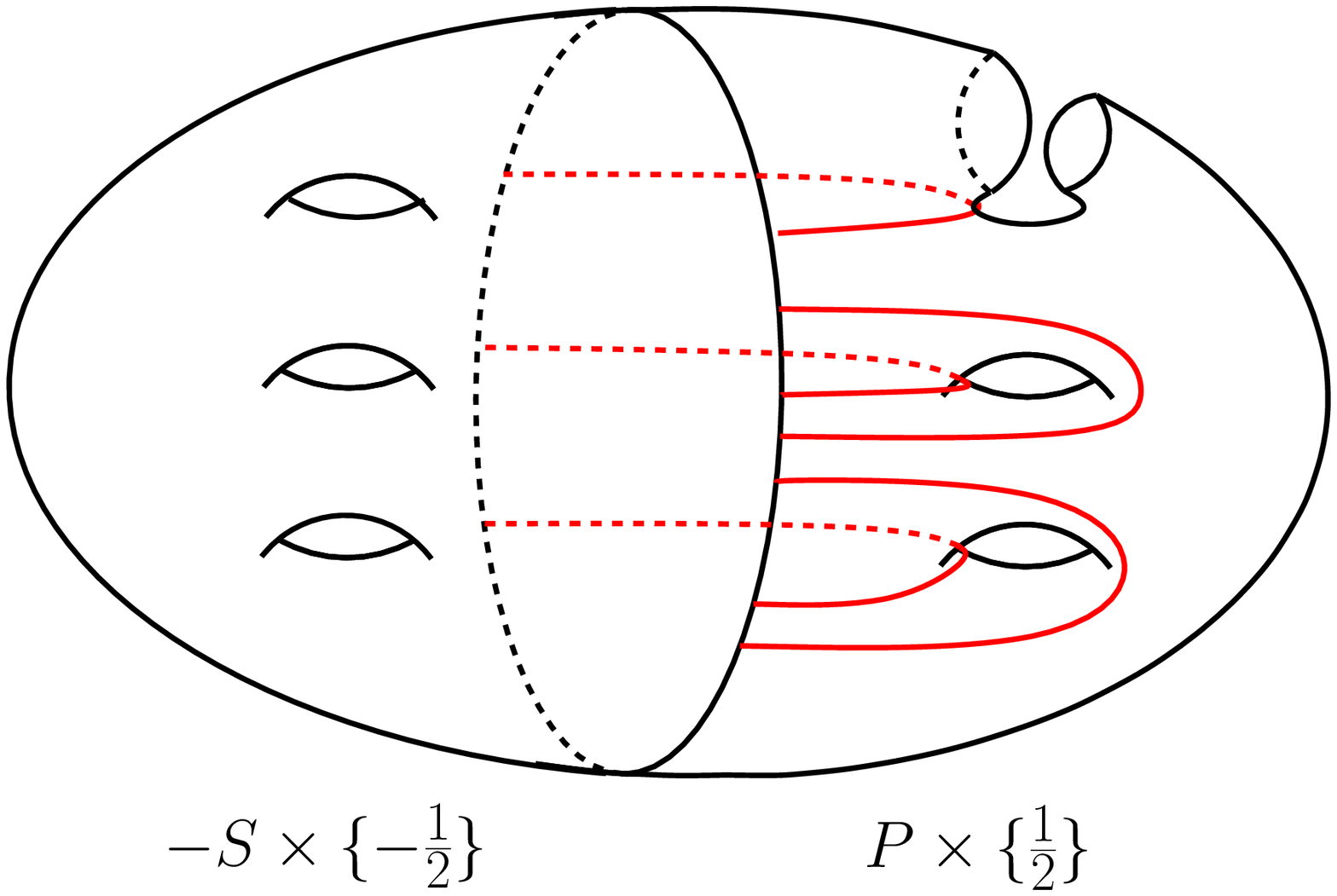}\\[0.2cm]
\mbox{Heegaard diagram for $Y$} & \mbox{Heegaard diagram for $Y-N(L)$}
\end{array}$
\caption{Heegaard diagrams corresponding to the (partial) open books.}
\label{fig:heegaardiagram}
\end{figure}

Attaching a bypass along $c$ changes the partial open book to
$(S^{\prime},P^{\prime},h^{\prime})$, where (with the notations described in
Subsection~\ref{subsec:bypass}) we have $S^{\prime}=S\cup(\textrm{1--handle})$
and $P^\prime=P\cup \nu (a_+)$. Note that $a_+$ represents half of the
meridian on $(\partial(\nu(L)))_+\subset S$, thus we can orient it.  The
1--handle is attached to $S$ along $\partial S$ in the neighbourhood of the
head of $a_+$ so that both of its feet are in the positive direction away from
the head of $a_+$ with respect to the orientation of $\partial S$, cf.
Figure~\ref{fig:heegaardsurface3}.  The monodromy remains the same restricted
to $P$ (i.e.\ $h^{\prime}|_P=h$) and as it was observed in Section
\ref{subsec:bypass}, $h^{\prime}(a_+)=a_-$ and $a_-$ splits as the core of the
1--handle and as $a_-\cap S$ which is isotopic to $c_-$. Note that $c_-$ is a
half--meridian of the knot $L$, thus the image of it on $S\times\{-\frac12\}$
is isotopic to $g(a_1)$.  Now we are ready to describe the Heegaard diagram
$(-(\Sigma^{\prime},\{\alpha,\alpha_2^{\prime}\dots,\alpha_k^{\prime}\},
\{\beta,\beta_2^\prime,\dots,\beta_k^\prime\}))$ obtained from the partial
open book $(S^\prime,P^\prime,h^\prime)$ in the usual manner. The Heegaard
surface $-\Sigma^\prime$ is equal to $P^\prime \times\{\frac12\} \cup
-S^\prime\times\{-\frac12\}$, and the curves
$\beta^\prime=b_+\times\{\frac12\}\cup b_+\times\{-\frac12\}$ and
$\alpha^\prime=a_+\times\{\frac12\} \cup a_-\times\{-\frac12\}$, where $b_+$
is the usual perturbation of $a_+$ on $P^\prime$.  $\Sigma^\prime$ is obtained
by gluing two surfaces together, each of which is diffeomorphic to
$S-\nu(\textrm{point})$. Indeed, the hole on the $S^\prime$--side comes from
the 1--handle attachment. $P^\prime$ is just a union of the 1--handles of $S$,
thus the missing 2--handle gives us the other hole. This surface
$\Sigma^\prime$ is thus diffeomorphic to $\Sigma-\nu(z)-\nu(w)$, where we
think of $\nu(z)$ being deleted from the $S^\prime$-- and $\nu({w})$ from the
$P^\prime$--side. Under this identification $b_+$ (and thus $a_+$) is isotopic
to $b_1$ on $P^\prime$, hence $\beta^\prime=b_+\times\{\frac12\}\cup
b_+\times\{-\frac12\}$ and $\beta_1$ are isotopic on $\Sigma^\prime$. Recall
that $h^\prime(a_+)$ on $S^\prime\times \{-\frac12\}$ was isotopic to the
union of $g(a_1)$ and the core of the 1--handle. So $\alpha^\prime$ is
isotopic to $\alpha_1$ on $\Sigma-\nu(z)$. The core part of $h^\prime(a_+)$
makes $\alpha^\prime$ and $\beta^\prime$ to go around the hole $\nu(w)$ from
different sides, thus $\alpha^\prime$ is isotopic to $\alpha_1$ on
$\Sigma^\prime$.
%The
%1--handle is attached to $\partial S$ close to the head of $a_+$ in the
%direction of the boundary orientation of $S$.  The monodromy $h^{\prime}$
%remains the same on $P$, and $h^\prime(a_+)=a_-$ splits as the core of the
%1--handle in $S^\prime$ and $a_-\cap S=h(c_-)$.
%
%Consider the surface $-\Sigma^\prime=P^\prime \times\{\frac12\} \cup
%-S^\prime\times\{-\frac12\}$. This is a surface with two boundary components,
%one is newly introduced and is on the $(-S^\prime)$--side, while the other is
%the boundary of the ``missing 0--handle on the $P^\prime$--side'', cf.
%Figure~\ref{fig:heegaardsurface3}. The new curves are
%$\beta^\prime=b_+\times\{\frac12\}\cup b_+\times\{-\frac12\}$ and
%$\alpha^\prime=a_+\times\{\frac12\} \cup a_-\times\{-\frac12\}$, where $b_+$
%is the usual perturbation of $a_+$ on $P^\prime$.  The curves $\alpha^\prime$
%and $\beta^\prime$ go around the hole on the $(-S^\prime)$--side from
%different directions.
In conclusion, the Heegaard diagram $(-\Sigma^\prime,
\{\alpha^\prime,\alpha_2,\dots,\alpha_k\},
\{\beta^\prime,\beta_2,\dots,\beta_k\})$ is isotopic to $(-(\Sigma-\nu
(z\cup w)),\{\alpha_1,\dots,\alpha_k\},\{\beta_1,\dots,\beta_k\})$.  The
contact invariant $\EH(L)$ is mapped to the contact invariant $\EH(Y-\nu
(L),-m\cup m)$ under the map induced by the basic slice, and thus it
represents the Legendrian invariant in $\CFKh(-\Sigma,\alphas,\betas,z,w)$,
which proves the statement.
\begin{figure}[ht]
\centering
\includegraphics[scale=0.4]{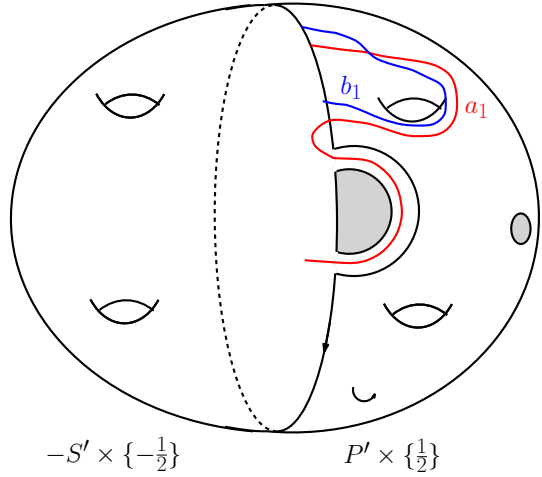}
\caption{Heegaard diagram corresponding to $(S^\prime,P^\prime,h^\prime)$.}
\label{fig:heegaardsurface3}
\end{figure}
\end{proof}
\begin{proof}[Proof of Theorem~\ref{t:main}]
With the identifications above, the proof of Theorem~\ref{t:main} is
now complete.
\end{proof}
\noindent Next we turn to the proof of the remaining statements
described in Section~\ref{sec:intro}.

\begin{proof}[Proof of Theorem \ref{thm:stabilization}]

  Take a standard contact neighbourhood $\nu (L)$ of $L$ and stabilize $L$
  inside it. Then $L^\pm$ has a standard contact neighbourhood $\nu
  (L^\pm)\subset \nu (L)$. As it is explained in \cite{EH}, the contact
  manifold $(\nu (L)-\nu (L^\pm),\xi|_{\nu (L)-\nu (L^\pm)})$ is a basic
  slice, i.e., $Y-\nu (L^{\pm})$ is obtained from $Y-\nu (L)$ by a bypass
  attachment. We can view $Y-\nu (L)$ as the result of a bypass attachment to
  the boundary of $Y-\nu (L^{\pm })$ \emph{from the back}. As usual, the two
  basic slices with the above boundary conditions have opposite relative Euler
  classes. To figure out which one corresponds to the positive and which one
  to the negative stabilization we first examine a model case. (For a related
  discussion see \cite{EH}.) Suppose that $\tb(L)<0$ and take a Seifert
  surface $S$ for $L$, giving rise to the Seifert surface $S^p$ (resp. $S^m$)
  for $L^+$ (resp.\ $L^-$). These surfaces are oriented such that their
  boundary orientations give the orientations for the knot. By $\tb (L)<0$ we
  can assume that $S$ is in convex position.  We have $\tb (L^\pm)=\tb(L)-1$,
  thus the dividing curve hits the boundary of the Seifert surface $S$ in
  $2|\tb(L)-1|$ points. In the collar neighbourhood of the boundary
  (diffeomorphic to $S^1\times I$), the dividing curves of $S$ are the line
  segments $k\frac{2\pi}{2|\tb(L)|}\times I$ where $0\leq k < 2|\tb(L)|$.
  Once again, by the negativity of $\tb (L)$ the bypass attachment corresponds
  to the gluing of an annulus to the boundary of $S$ with dividing curves
  $k\frac{2\pi}{2|\tb(L)|}\times I$ ($0\leq k < 2|\tb(L)|$) and a boundary
  parallel curve that is disjoint from the others. This boundary parallel
  curve bounds a domain, cf.  Figure~\ref{fig:stabilization1}.  The rotation
  numbers are $\rot (L^\pm)=\rot (L)\pm 1$, thus by the formula $\rot
  (S)=\chi(S_+)-\chi(S_-)$ we get that the extra domain on $S^p$ (on $S^m$,
  resp.) is in the positive (resp.  negative) region.  Using edge rounding we
  get that the attaching curve corresponding to the positive (resp. negative)
  stabilization must end in the positive (resp. negative) region with respect
  to the orientation of the knot.  The left hand side of Figure
  \ref{fig:stabilization2} depicts the arc $p$ (and $n$, resp.) along which
  the bypass has to be attached (from the back) to obtain $Y-\nu (L)$.

\begin{figure}[ht]
\centering \includegraphics[scale=0.35]{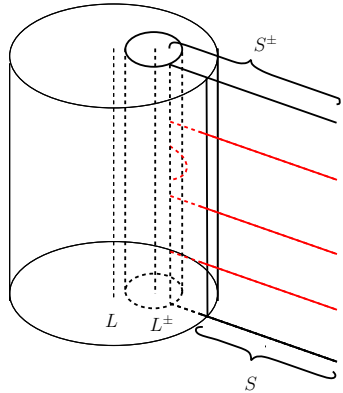}
\caption{Neighbourhood of a Legendrian knot and its stabilization.}
\label{fig:stabilization1}
\end{figure}

Both the stabilization and the bypass attachment are local operations, thus
the above described phenomenon remains true for any Legendrian knot (without
the assumption $\tb (L)<0$).  The arcs $p$ and $n$ have the same slope, but
they end in regions of different sign.  Consider the middle diagram of
Figure~\ref{fig:stabilization2} for the general picture for $T^2$, trivialized
by the meridian $m$ and the Thurston--Bennequin framing $l$.

By Theorem~\ref{thm:gluing} the map corresponding
to the bypass attachment maps $\EH(L)$ to
$\EH(L^\pm)$. To get $\LegInva (L^\pm)$
we need to attach another bypass, so that the new dividing curves are
meridians, hence this second bypass is attached along the arc $c$.

In the case of positive stabilization, the manifold $(Y-\nu (L^+),(\xi|_{Y-\nu
  (L^+)})^c)=(Y-\nu (L),(\xi|_{Y-\nu (L)})^{p^{-1}c})$ is overtwisted.
Indeed, performing the positive stabilization first one can indicate both
bypasses in one picture, one attached from the back: $p^{-1}$ drawn by dashed
line on Figure~\ref{fig:stabilization2} and $c$ from the front. These curves
are parallel, thus the corresponding bypasses ('half overtwisted disks') form
an overtwisted disk in $(Y-\nu (L),(\xi|_{Y-\nu (L)})^{p^{-1}c})$. It is known
that the sutured invariant of an overtwisted structure vanishes
\cite[Corollary 4.3.]{HKM2}, therefore so does $\LegInva (L^+)$.

In the case of negative stabilization, the contact structure $(T^2\times
I,\xi^{n^{-1}c})$ is universally tight.  This can be seen by first passing to
$\partial(Y-\nu (L))$ (cf. the right hand side of
Figure~\ref{fig:stabilization2}) and then noting that the two bypasses
attached there are of the same sign, so they do not induce an overtwisted
disk. The union of the two basic slices is minimally twisting, and in this
case the range of slopes is $[0,\infty]=[0,1]\cup[1,\infty]$.  Therefore the
result is still a basic slice, thus the composition of the two bypass
attachments along $n$ and $c$ is equivalent to a single bypass attachment
along $c$. This immediately implies $\LegInva (L)=\LegInva (L^-)$, concluding
the proof.
\begin{figure}[ht]
\centering \includegraphics[scale=0.7]{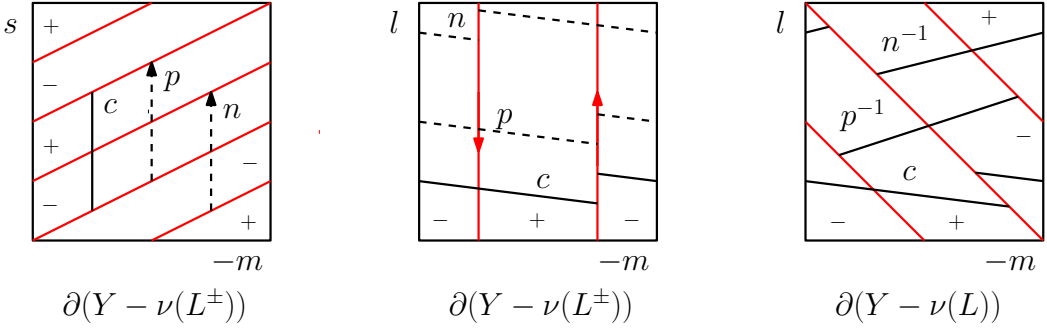}
\caption{Attaching curves for the bypasses corresponding to the
  stabilizations. The dashed line indicates that the bypass is attached from
  the back. On the left--hand picture $s$ denotes the Seifert framing of the
  knot, while on the two right--hand pictures $l$ is given by the contact
  framing of the Legendrian knot.}
\label{fig:stabilization2}
\end{figure}
\end{proof}

\noindent Next we turn to the proof of the statement concerning the vanishing
of the Legendrian invariant in the presence of Giroux torsion. We start
by recalling Giroux torsion.

\begin{definition} {\rm The contact structure $\xi_n$ on
    $T^2\times [0,1]=\mathbb{R}/\mathbb{Z}\times\mathbb{R}/\mathbb{Z}\times
    [0,1]=\{(x,y,z))\}$ is defined by $\xi_n=\ker(\cos(2\pi
    nz)\mathit{d}x-\sin(2\pi nz)\mathit{d}y)$.  A (not necessarily closed)
    contact 3--manifold $(Y,\xi)$ has \emph{Giroux torsion} $\tau(Y,\xi)\ge n$
    if it contains an embedded submanifold $T^2\times I$ with the property
    that $(T^2\times I , \xi \vert _{T^2\times I})$ is contactomorphic to
    $(T^2\times [0,1],\xi_n)$.}
\end{definition}

\begin{proof}[Proof of Corollary~\ref{cor:girouxtorsion}]
The proof is a simple adaptation of the proof for the closed case
given by Ghiggini, Honda, and Van Horn-Morris \cite{GHV}.  As $(Y-\nu
(L),\xi_{Y-\nu (L)})$ has positive Giroux torsion, there is a
submanifold $T^2\times I$, such that $\xi|_{T^2\times I}=\xi_n$ for
some $n>0$. It was shown in \cite{GHV} that $\EH(T^2\times
I,\xi_n)=0$.

The application of Theorem~\ref{thm:gluing} for the contact
3--manifold pair $(Y-\nu (L), T^2\times [0,1])$ provides a map
\[
\SFH(-(T^2\times I),-\Gamma_{\partial (T^2\times I)})
\to\SFH(-(Y-\nu (L)),-\Gamma_{Y- \nu (L)}))
\]
mapping the contact element $\EH(T^2\times I,\xi_n)=0$ to the contact
element $\EH (L) =\EH(Y-\nu (L),\xi|_{Y-\nu (L)})$.  This implies that
$\EH (L)=0$, hence in the light of Theorem~\ref{t:main} we get that
$\LegInva (L)=0$, concluding the proof.
\end{proof}

\begin{proof}[Proof of Theorem \ref{thm:surgery}]

  As in the proof of Theorem \ref{t:main}, we attach a bypass along the arc
  $e$ of Figure~\ref{fig:surgery} and change the dividing curve
  on the torus boundary to $\Gamma_{\partial(Y-\nu(L))}^{e}$ of slope $-1$.
  There are two choices for such arcs, but again the orientation of $L$
  assigns the one depicted on Figure \ref{fig:surgery}.

\begin{figure}[ht]
\centering \includegraphics[scale=0.7]{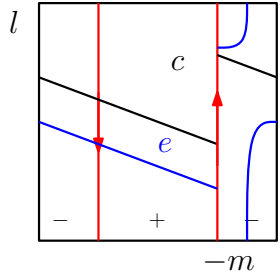}
\caption{Attaching curves for the bypasses on $\partial(Y-\nu(L))$ to obtain dividing
  curves of slope $1$.}
\label{fig:surgery}
\end{figure}

 This bypass attachment gives rise to a map
\[
\Phi^{e}\colon \SFH(-(Y-\nu(L)),-\Gamma_{\partial(Y-\nu(L))})\to
\SFH(-(Y-\nu(L)),-\Gamma_{\partial(Y-\nu(L))}^{e}).
\]
By filling the boundary with a solid torus, the latter homology is identified
with $\HFKh(-Y_{-1},L^\prime)$. Denote the composition of the above maps by
\[
G\colon \SFH(-(Y-\nu(L)),-\Gamma_{\partial(Y-\nu(L))})\to
\HFKh(-Y_{-1},L^\prime).
\]
We claim that the homomorphism $G$ maps $\EH(L)$ to $\LegInva(L^\prime)$.
Indeed, consider an open book $(S,h)$ adapted to $(Y, \xi , L )$. The same
open book is adapted to $(Y_{-1}, \xi _{-1}, L ')$, with the only difference
in the monodromy: the monodromy $h'$ for the latter triple is multiplied by a
right--handed Dehn twist along $L$, cf. \cite[page~133]{OzS}. Using the
notations introduced in Section \ref{sec:invariant}, the map $G$ corresponds
to changing the partial open book $(S,P=S-\nu _S(L),h|_{P})$ to
$(S^\prime,P^\prime,h^{\prime\prime})$ corresponding to the bypass attachment.
The image of the half--meridian $a_+$ under $h^{\prime\prime}$ is $h(a^+)$
multiplied by a right--handed Dehn twist along $L$. Therefore
$G(\EH(L))=\LegInva(L^\prime)$.

After attaching the bypass along $e$, we can apply another bypass attachment
along $c$ of Figure~\ref{fig:surgery} to obtain the meridian as
dividing curve.  We have already seen in the proof of
Theorem~\ref{thm:stabilization} that the composition of these two bypasses is
a basic slice, thus we have the commutative diagram

\xymatrix{
*\txt<14pc>{%
$\SFH(-(Y-\nu(L)),-\Gamma_{\partial(Y-\nu(L))})$} \ar@{->}[dr]
&{}\save[]+<2cm,0cm>*\txt<15pc>{%
$\SFH(-(Y-\nu(L)),-\Gamma^e_{\partial(Y-\nu(L))})=\HFKh(-Y_{-1},L^\prime)$}
\ar@{<-}[l] \ar[d] \restore \\
& *\txt<13pc>{%
 $\SFH(-(Y-\nu(L)),-m\cup m)=\HFKh(-Y,L)$}   &  &
}

%\[
%\xymatrix{
%\SFH(-(Y-\nu(L)),-\Gamma_{\partial(Y-\nu(L))}) \ar@/^2pc/[rr] \ar[r]
%& \SFH(-(Y-\nu(L)),-\Gamma_{\partial(Y-\nu(L))}^{e}) \ar[r] \ar@{=}[d]
%& \SFH(-(Y-\nu(L)),-m\cup m) \ar@{=}[d]  \\
%& \HFKh(-Y_{-1},L^\prime) & \HFKh(-Y,L) \\
% }
%\]

The maps in the above triangle map the contact invariants as
%\[
%\xymatrix{
%\EH(L) \ar@/^1pc/[rr] \ar[r]
%& \LegInva(L^\prime)\ar[r]
%& \LegInva(L)\neq 0,   \\
% }
%\]

\[\xymatrix{\EH(L)\ar[rr]\ar@{->}[dr] & &\LegInva(L^\prime)\ar[dl]\\
 &\LegInva(L)\neq 0 &  }
\]

therefore $\LegInva(L^\prime)$ does not vanish, concluding the proof.
\end{proof}

\bibliographystyle{plain}

\end{document}

%% file: bas.tex
\newtheorem {theorem}{Theorem}[section]
\newtheorem {definition}[theorem]{Definition}

\newtheorem {corollary}[theorem]{Corollary}

\newtheorem {remark}[theorem]{Remark}

\def\EH{\mathrm{EH}}

\newcommand\SFH{\mathrm{SFH}}
\newcommand\SFC{\mathrm{SFC}}
\newcommand\CFKh{\widehat{\mathrm{CFK}}}
\newcommand\HFKh{\widehat{\mathrm{HFK}}}
\newcommand\HFKa{\widehat{\mathrm{HFK}}}
\newcommand\CFKa{\widehat{\mathrm{CFK}}}
\newcommand\HFa{\widehat{\mathrm{HF}}}
\newcommand\CFa{\widehat{\mathrm{CF}}}

\newcommand{\x}{\mathbf{x}}
\newcommand{\y}{\mathbf{y}}

\newcommand\betas{\mbox{\boldmath$\beta$}}
\newcommand\alphas{\mbox{\boldmath$\alpha$}}

\newcommand\Sigmah{\widetilde{\Sigma}}

\newcommand\LegInva{\widehat{\mathcal{L}}}
\newcommand\tb{{\rm {tb}}}
\newcommand\rot{{\rm {rot}}}
\newcommand\bfz{\Bbb {Z}}
\newcommand\bfr{\Bbb {R}}
\newcommand\Ta{{\Bbb {T}}_{\alpha}}
\newcommand\Tb{{\Bbb {T}}_{\beta}} 